\documentclass[]{interact}

\usepackage{epstopdf}
\usepackage[caption=false]{subfig}

\usepackage[numbers,sort&compress]{natbib}
\bibpunct[, ]{[}{]}{,}{n}{,}{,}
\makeatletter
\def\NAT@def@citea{\def\@citea{\NAT@separator}}
\makeatother

\theoremstyle{plain}
\newtheorem{theorem}{Theorem}
\newtheorem{lemma}[theorem]{Lemma}
\newtheorem{corollary}[theorem]{Corollary}
\newtheorem{proposition}[theorem]{Proposition}

\theoremstyle{definition}

\theoremstyle{remark}

\def \la {\lambda}
\def \al {\alpha}
\def \be {\beta}

\begin{document}


\title{On the extreme zeros of Jacobi polynomials}

\author{
\name{Geno Nikolov\thanks{CONTACT Geno Nikolov. Email:
geno@fmi.uni-sofia.bg}} \affil{Faculty of Mathematics and
Informatics, Sofia University St. Kliment Ohridski, 5 James
Bourchier blvd., 1164 Sofia, Bulgaria}}

\maketitle

\begin{abstract}
By applying the Euler--Rayleigh methods to a specific representation
of the Jacobi polynomials as hypergeometric functions, we obtain new
bounds for their largest zeros. In particular, we derive upper and
lower bound for $1-x_{nn}^2(\la)$, with $x_{nn}(\la)$ being the
largest zero of the $n$-th ultraspherical polynomial $P_n^{(\la)}$.
For every fixed $\la>-1/2$, the limit of the ratio of our upper and
lower bounds for $1-x_{nn}^2(\la)$ does not exceed $1.6$. This paper
is a continuation of \cite{GN2019}.
\end{abstract}

\begin{keywords}
Jacobi polynomials; Gegenbauer polynomials; Laguerre polynomials;
Euler-Rayleigh method
\end{keywords}

\begin{amscode}
33C45; 42C05
\end{amscode}

\section{Introduction and statement of the results}
The extreme zeros of the classical orthogonal polynomials of Jacobi,
Laguerre and Hermite have been a  subject of intensive study. We
refer to Szeg\H{o}'s monograph \cite{GS1975} for earlier results,
and to \cite{IsmLi(1992),IsmMul(1995),Kras(2003),{ADGR(2004)},
GupMul(2007),DimNik(2010),ADGR(2012),DriJor(2012),DriJor(2013),
GN2019,NU(2019)} for some recent developments.

\label{sec:1} Throughout this paper we use the notation
$$
x_{1n}(\al,\be)<x_{2n}(\al,\be)<\cdots<x_{nn}(\al,\be)
$$
for the zeros of the $n$-th Jacobi polynomial $P_n^{(\al,\be)}$,
$\al,\be>-1$, and the zeros of the $n$-th Gegenbauer polynomial
$P_n^{(\la)}$, $\la>-1/2$ are denoted by
$$
x_{1n}(\la)<x_{2n}(\la)<\cdots<x_{nn}(\la)\,.
$$

In the recent paper \cite{GN2019} we applied the Euler--Rayleigh
method to the Jacobi and, in particular, the Gegenbauer polynomials,
represented as hypergeometric functions, to derive new bounds for
their extreme zeros. Below we state some of the bounds obtained in
\cite{GN2019}, which improve upon some results of Driver and Jordaan
\cite{DriJor(2012)}.
\medskip

\noindent \textbf{Theorem A. (\cite[Theorem 1.4]{GN2019})} \emph{For
every $n\geq 3$ and $\al,\,\be
>-1$, the largest zero $\,x_{nn}(\al,\be)\,$ of the Jacobi
polynomial $P_n^{(\al,\be)}$ satisfies
$$
1-x_{nn}(\al,\be)<\frac{2(\al+1)(\al+3)}
{(n+\al+1)(n+\al+\be+1)\big[2-\frac{(\al+1)(2n+\be-1)}
{(n+\al+1)(n+\al+\be+1)-(\al+1)(\al+2)}\big]}.
$$}
\medskip

\noindent \textbf{Corollary A. (\cite[Corollary 1.6]{GN2019})}
\emph{For every $n\geq 3$ and $\la>-1/2$, the largest zero
$\,x_{nn}(\la)\,$ of the Gegenbauer polynomial $P_n^{(\la)}$
satisfies
$$
1-x_{nn}(\la)<\frac{(2\la+1)(2\la+5)}
{(n+2\la)(2n+2\la+1)\big[2-\frac{(2\la+1)(4n+2\la-3)}
{2(n+2\la)(2n+2\la+1)-(2\la+1)(2\la+3)}\big]}.
$$}

\noindent \textbf{Theorem B. (\cite[Theorem 1.1]{GN2019})} \emph{For
every $n\geq 3$ and $\la>-1/2$, the largest zero $\,x_{nn}(\la)\,$
of the Gegenbauer polynomial $P_n^{(\la)}$ satisfies
$$
1-x_{nn}^2(\la)<\frac{(2\la+1)(2\la+5)}{2n(n+2\la)+2\la+1
+\frac{2(\la+1)(2\la+1)^2(2\la+3)}{n(n+2\la)+2(2\la+1)(2\la+3)}}\,.
$$}
\medskip

The above results provide lower bounds for the largest zeros of the
Jacobi and Gegenbauer polynomials. It is instructive to compare
Theorem~B with the following upper bound for the largest zeros of
the Gegenbauer polynomials:
\medskip

\noindent \textbf{Theorem C. (\cite[Lemma 3.5]{GN2005})} \emph{For
every  $\la>-1/2$, the largest zero $\,x_{nn}(\la)\,$ of the
Gegenbauer polynomial $P_n^{(\la)}$ satisfies
$$
1-x_{nn}^2(\la)>\frac{(2\la+1)(2\la+9)}{4n(n+2\la)+(2\la+1)(
2\la+5)}\,.
$$}
\medskip

We observe that, for any fixed $\la>-1/2$ and large $n$, the ratio
of the upper and the lower bound for $1-x_{nn}^2(\la)$, given by
Theorems B and C, does not exceed $2$. With Corollary~\ref{c1} below
this ratio is reduced to $1.6$.

In the present paper we apply the Euler--Rayleigh method to the
Jacobi polynomial $P_n^{(\al,\be)}$, represented as a hypergeometric
function, to obtain further bounds for the largest zeros of the
Jacobi and Gegenbauer polynomials. As at some points the
calculations become unwieldy, we have used the assistance of
\emph{Wolfram Mathematica}.

The following is the main result in this paper.

\begin{theorem} \label{t1}
For every $n\geq 4$ and $\al,\,\be
>-1$, the largest zero $\,x_{nn}(\al,\be)\,$ of the Jacobi
polynomial $P_n^{(\al,\be)}$ satisfies
\begin{equation}\label{e1}
1-x_{nn}(\al,\be)<\frac{4(\al+1)(\al+2)(\al+4)}
{(5\al+11)\big[n(n+\al+\be+1)+\frac{1}{3}(\al+1)(\be+1)\big]} \,.
\end{equation}
Moreover, if either $\,n\geq\max\{4,\al+\be+3\}\,$ or $\,\be\leq
4\al+7$, then
\begin{equation}\label{e2}
1-x_{nn}(\al,\be)<\frac{4(\al+1)(\al+2)(\al+4)}
{(5\al+11)\big[n(n+\al+\be+1)+\frac{1}{2}(\al+1)(\be+1)\big]}\,.
\end{equation}
\end{theorem}
Since $P_n^{(\al,\be)}(x)=(-1)^{n}P_n^{(\be,\al)}(-x)$,
Theorem~\ref{t1} can be equivalently formulated as

\begin{theorem} \label{t2}
For every $n\geq 4$ and $\al,\,\be
>-1$, the smallest zero $\,x_{1n}(\al,\be)\,$ of the Jacobi
polynomial $P_n^{(\al,\be)}$ satisfies
$$
1+x_{1n}(\al,\be)<\frac{4(\be+1)(\be+2)(\be+4)}
{(5\be+11)\big[n(n+\al+\be+1)+\frac{1}{3}(\al+1)(\be+1)\big]} \,.
$$
Moreover, if either $\,n\geq\max\{4,\al+\be+3\}\,$ or $\,\al\leq
4\be+7$, then
$$
1+x_{1n}(\al,\be)<\frac{4(\be+1)(\be+2)(\be+4)}
{(5\be+11)\big[n(n+\al+\be+1)+\frac{1}{2}(\al+1)(\be+1)\big]}\,.
$$
\end{theorem}
The assumption $\be\leq 4\al+7$ is satisfied, in particular, when
$\be=\al>-1$. Therefore, as a consequence of Theorem~\ref{t1}, we
obtain a bound for the largest zero of the ultraspherical polynomial
$P_n^{(\la)}=c\,P_n^{(\al,\al)}$, $\,\al=\la-\frac{1}{2}$.

\begin{theorem} \label{t3}
For every $n\geq 4$ and $\la>-1/2$, the largest zero
$\,x_{nn}(\la)\,$ of the Gegenbauer polynomial $P_n^{(\la)}$
satisfies
\begin{equation}\label{e3}
1-x_{nn}(\la)<\frac{(2\la+1)(2\la+3)(2\la+7)}
{(10\la+17)\big[n(n+2\la)+\frac{1}{8}(2\la+1)^2\big]}\,.
\end{equation}
\end{theorem}
Theorem~\ref{t3} and $\,1-x_{nn}^2(\la)<2(1-x_{nn}(\la))\,$ imply
immediately the following:

\begin{corollary}\label{c1}
For every $n\geq 4$ and $\la>-1/2$, the largest zero
$\,x_{nn}(\la)\,$ of the Gegenbauer polynomial $P_n^{(\la)}$
satisfies
\begin{equation}\label{e4}
1-x_{nn}^2(\la)<\frac{2(2\la+1)(2\la+3)(2\la+7)}
{(10\la+17)\big[n(n+2\la)+\frac{1}{8}(2\la+1)^2\big]}\,.
\end{equation}
\end{corollary}

Usually, the comparison of the various bounds for the extreme zeros
of the classical orthogonal polynomials is not an easy task due to
the parameters involved. At least for large $n$, the bounds provided
by Theorem~\ref{t1}, Theorem~\ref{t3} and Corollary~\ref{c1} are
sharper than those in Theorem~A, Corollary~A and Theorem~B,
respectively. In fact, the actual bounds obtained with the approach
here are slightly sharper but are given by rather complicated
expressions; in particular, by a limit passage we reproduce a result
of Gupta and Muldoon from \cite{GupMul(2007)} concerning the
smallest zero of the Laguerre polynomial. These and some other
observations are given in Section~\ref{s4} of the paper.

The rest of the paper is organized as follows. In Section~\ref{s2}
we present the necessary facts about the Euler--Rayleigh method and
the Newton identities. The proof of Theorem~\ref{t1} is given in
Section~\ref{s3.1}. For the reader's convenience, in
Section~\ref{s3.2} we include a short proof of Theorem~C.
\section{The Euler--Rayleigh method}\label{s2}
As was already mentioned, the proof of our results exploits the
so-called Euler--Rayleigh method (see \cite{IsmMul(1995)}). Here
(and also in \cite{GN2019}) the Euler--Rayleigh method is applied to
real-root polynomials, and for the reader's convenience we provide
some details from \cite{GN2019}.

Let $P$ be a monic polynomial of degree $n$ with zeros $(x_i)_1^n$,
\begin{equation}\label{e5}
P(x)=x^n-b_{1}\,x^{n-1}+b_{2}\,x^{n-2}-\cdots+(-1)^{n}b_n
=\prod_{i=1}^{n}(x-x_i)\,.
\end{equation}
For $k\in \mathbb{N}_0$, the power sums
$$
p_k=p_k(P):=\sum_{i=1}^{n}x_{i}^k,\qquad p_0=n={\rm deg\,}P,
$$
and the coefficients $(b_i)_1^n$ of $P$ are connected by the Newton
identities (cf. \cite{Waer(1949)})
$$
p_r+\sum_{i=1}^{\min\{r-1,n\}}(-1)^{i}p_{r-i}\,b_i+(-1)^r
r\,b_r=0\,.
$$
From Newton's identities one easily obtains:

\begin{lemma}\label{l1}
Assuming $n\geq r$, the following formulae hold for $p_r$, $1\leq
r\leq 4$:
\begin{align*}
  p_1(P) & = b_1\,; \\
  p_2(P) & = b_1^2-2b_2\,; \\
  p_3(P) & = b_1^3-3b_1b_2+3b_3\,; \\
  p_4(P) & = b_1^4-4b_1^2b_2+2b_2^2+4b_1b_3-4b_4\,.
\end{align*}
\end{lemma}

Let us set
$$
  \ell_k(P) := \frac{p_k(P)}{p_{k-1}(P)}\,, \qquad
  u_k(P) := \big[p_k(P)\big]^{1/k}\,, \qquad k\in \mathbb{N}\,.
$$
The following statement is Proposition 2.2 in \cite{GN2019}; it is a
slight modification of Lemma~3.2 in \cite{IsmMul(1995)}.

\begin{proposition} \label{p1}
Let $\,P$ be as in \eqref{e5} with positive zeros $\,x_1<
x_2<\cdots< x_n\,$. Then the largest zero $\,x_n\,$ of $\,P\,$
satisfies the inequalities
$$
  \ell_k(P) < x_n < u_k(P)\,, \qquad k\in \mathbb{N}\,.
$$
Moreover, $\,\{\ell_k(P)\}_{k=1}^{\infty}\,$ is monotonically
increasing, $\,\{u_k(P)\}_{k=1}^{\infty}\,$ is monotonically
decreasing, and
$$
\lim_{k\to\infty} \ell_k(P) = \lim_{k\to\infty} u_k(P) = x_n\,.
$$
\end{proposition}
\section{Proof of the results}
\subsection{Proof of Theorem~\ref{t1}}
\label{s3.1} The starting point for the proof of Theorem~\ref{t1} is
the following representation of $P_n^{(\al,\be)}$ (cf. \cite[eqn.
(4.21.2)]{GS1975}:
\begin{equation}\label{e6}
P_n^{(\al,\be)}(x)=\frac{(\al+1)_n}{n!}\,
_2F_1\Big(-n,n+\al+\be+1;\al+1;\frac{1-x}{2}\Big)
\end{equation}
(for the proof of Theorem~A we have used another representation of
$P_n^{(\al,\be)}$ as a hypergeometric function, namely, \cite[eqn.
(4.3.2)]{GS1975}). Here we use Szeg\H{o}'s notation for the
hypergeometric $_2F_1$ function,
$$
F(a,b;c;z)=1+\sum_{k=1}^{\infty}\frac{(a)_k}{k!}\,\frac{(b)_k}{(c)_k}\,z^k\,,
\qquad (a)_k:=a(a+1)\cdots(a+k-1)\,.
$$

It follows from \eqref{e6} that the monic polynomial
$$
P(z)=z^n+\sum_{i=1}^{n}(-1)^{i}b_i\,z^{n-i}
$$
with coefficients
\begin{equation}\label{e7}
b_i=b_i(P)={n\choose i}\frac{(n+\al+\be+1)_i}{(\al+1)_i}\,,\qquad
i=1,\ldots,n,
\end{equation}
has $n$ positive zeros $z_1<z_2<\cdots<z_n$, connected with the
zeros of $P_{n}^{(\al,\be)}$ by the relation
$$
z_i=\frac{2}{1-x_{in}(\al,\be)}\,,\qquad i=1,\ldots,n\,.
$$
According to Proposition~\ref{p1},
$p_{k+1}(P)/p_{k}(P)<z_n<\big[p_k(P)\big]^{1/k}$, $k\in \mathbb{N}$,
and consequently
\begin{equation}\label{e8}
\frac{2}{\big[p_k(P)\big]^{1/k}}<1-x_{nn}(\al,\be)<
\frac{2p_{k}(P)}{p_{k+1}(P)}\,,\qquad k\in \mathbb{N}\,.
\end{equation}
At this point, we find it suitable to substitute
\begin{eqnarray*}
&&a:=\al+1\,,\ b:=\be+1\,,\\
&&t:=n(n+\al+\be+1)\,,
\end{eqnarray*}
thus $\,a,\,b>0\,$ and $\,t=n(n+a+b-1)\,$. With this notation, the
first four coefficients $\,b_i(P)\,$ in \eqref{e7} are given by
\begin{eqnarray*}
&& b_1(P)=\frac{t}{a}\,,\quad
b_2(P)=\frac{t(t-a-b)}{2a(a+1)}\,,\quad
b_3(P)=\frac{t(t-a-b)\big[t-2(a+b+1)\big]}{6a(a+1)(a+2)}\,,\\
&& b_4(P)=\frac{t(t-a-b)\big[t-2(a+b+1)\big]\big[t-3(a+b+2)\big]}
{24a(a+1)(a+2)(a+3)}\,.
\end{eqnarray*}
Using Lemma~\ref{l1}, we find $p_1(P)=b_1(P)=t/a$,
\begin{equation}\label{e9}
p_2(P)=\frac{t\big[t+a(a+b)\big]}{a^2(a+1)}\,,
\end{equation}
\begin{equation}\label{e10}
\begin{split}
p_3(P)=&\frac{t\,q_2(t)}{a^3(a+1)(a+2)}\,,\\
q_2(t)=&2t^2+a(2a+3b)t+a^2(a+b)(a+b+1)\,,
\end{split}
\end{equation}
\begin{equation}\label{e11}
\begin{split}
p_4(P)=&\frac{t\,q_3(t)} {a^4(a+1)^2(a+2)(a+3)}\,,\\
q_3(t)=&(5a+6)t^3+2a(3a^2+5ab+4a+6b)t^2\\
&+a^2(3a^3+9a^2b+6ab^2+6a^2+15ab+7b^2+2a+4b)t\\
&+a^3(a+1)(a+b)(a+b+1)(a+b+2)\,.
\end{split}
\end{equation}

Theorem~\ref{t1} follows from the right-hand inequality in
\eqref{e8} with $k=3$. In order to show this, we observe that,
according to \eqref{e10} and \eqref{e11},
$$
1-x_{nn}(\al,\be)< \frac{2p_{3}(P)}{p_{4}(P)}
=\frac{2a(a+1)(a+3)q_2(t)}{q_3(t)}=\frac{4(\al+1)(\al+2)(\al+4)}
{\frac{2q_3(t)}{q_2(t)}}\,.
$$
Hence, to prove the first part of Theorem~\ref{t1}, it suffices to
show that if $a,\,b$ and $\,t\,$ are positive, then
$$
\frac{2q_3(t)}{q_2(t)}\geq (5a+6)\Big(t+\frac{ab}{3}\Big)
=(5\al+11)\Big[n(n+\al+\be+1)+\frac{1}{3}(\al+1)(\be+1)\Big]\,.
$$
With the help of \emph{Wolfram Mathematica} we find
$$
2q_3(t)-(5a+6)\Big(t+\frac{ab}{3}\Big)q_2(t)=\frac{a}{3}\,r_2(t)\,,
$$
where
\[
\begin{split}
r_2(t)=&(6a^2+5ab+12a+6b)t^2\\
&+a(3a^3+14a^2b+6ab^2+3a^2+27ab+6b^2-6a+6b)t\\
&+a^2(a+b)(a+b+1)(6a^2+ab+18a+12)\,.
\end{split}
\]
It is clear now that $r_2(t)>0$: the single negative summand in the
right-hand side, $-6a^2t$, is neutralized by $6a^2t^2$, since
$t=n(n+a+b-1)>n(n-1)\geq 12$ for $n\geq 4$. Consequently,
$$
\frac{2q_3(t)}{q_2(t)}\geq (5a+6)\Big(t+\frac{ab}{3}\Big)\,,
$$
and the first claim of Theorem~\ref{t1} is proved.

For the proof of the second claim of Theorem~\ref{t1} we need to
show that
\begin{equation}\label{e12}
2q_3(t)-(5a+6)\Big(t+\frac{ab}{2}\Big)q_2(t)\geq 0
\end{equation}
provided either $\,\be\le 4\al+7\,$ or $\,n\geq \max\{4,\al+\be+3\}=
\max\{4,a+b+1\}\,$.

With the assistance of \emph{Mathematica} we find
$$
2q_3(t)-(5a+6)\Big(t+\frac{ab}{2}\Big)q_2(t)=\frac{1}{2}a^2(a+2)\,s_2(a,b;t)\,,
$$
where
$$
s_2(a,b;t)=4t^2+(2a^2+6ab-b^2-2a+2b)t+a(a+b)(a+b+1)(4a-b+4)\,.
$$

Firstly, assume that $\,\be\leq 4\al+7\,$, which is equivalent to
$\,b\leq 4a+4$\,. Then obviously the constant term in the quadratic
$\,s_2(a,b;\cdot)\,$ is non-negative. We shall prove that the sum of
the other two terms is positive. Indeed, since for $n\geq 4$ we have
$$
t=n(n+a+b-1)\geq 4(a+b+3)>0\,,
$$
we need to show that $\,4t+2a^2+6ab-b^2-2a+2b>0$\,. This inequality
follows from
\[
\begin{split}
4t+2a^2+6ab-b^2-2a+2b&\geq
16(a+b+3)+2a^2+6ab-b^2-2a+2b\\
&=2a^2+b(6a+18-b)+14a+48\\
&\geq 2a^2+b(2a+14)+14a+48>0\,.
\end{split}
\]

Secondly, assume that $n\geq \max\{4,\al+\be+3\}= \max\{4,a+b+1\}$.
We observe that
\begin{equation}\label{e13}
t=n(n+a+b-1)\geq 2(a+b)(a+b+1)\,.
\end{equation}
Therefore,
\[
\begin{split}
4t^2\!+\!(2a^2\!+\!6ab\!-\!b^2\!-\!2a\!+\!2b)t&\geq
8(a\!+\!b)(a\!+\!b\!+\!1)t\!+\!(2a^2\!+\!6ab\!-\!b^2\!-\!2a\!+\!2b)t\\
&= (10a^2+22ab+7b^2+6a+10b)t>0 \,.
\end{split}
\]
Using this last inequality and applying \eqref{e13} once again, we
conclude that
\[
\begin{split}
s_2(a,b;t)&\geq (10a^2+22ab+7b^2+6a+10b)t+a(a+b)(a+b+1)(4a-b+4)\\
&\geq (a+b)(a+b+1)\Big[20a^2+44ab+14b^2+12a+20b+a(4a-b+4)\Big]\\
&=(a+b)(a+b+1)(24a^2+43ab+14b^2+16a+20b)>0\,.
\end{split}
\]
Thus, \eqref{e12} holds true in the case $n\geq
\max\{4,\al+\be+3\}$, which completes the proof of the second claim
of Theorem~\ref{t1}.
\subsection{Proof of Theorem~C}\label{s3.2}
The original proof of Theorem~C in \cite{GN2005} makes use of an
idea from \cite[Paragraph 6.2]{GS1975}, based on the following
observation of Laguerre: if $f$ is a real-valued polynomial of
degree $n$ having only real and distinct zeros, and $f(x_0)=0$, then
\begin{equation}\label{e14}
3(n-2)\big[f^{\prime\prime}(x_0)\big]^2
-4(n-1)f^{\prime}(x_0)f^{\prime\prime\prime}(x_0)\geq 0\,.
\end{equation}

In \cite{NU(2004)} Uluchev and the author proved a conjecture of
Foster and Krasikov \cite{FK(2002)}, stating that if $f$ is a
real-valued polynomial of degree $n$, then for every integer $m$
satisfying $0\leq 2m\leq n$ the following inequalities hold true:
$$
\sum_{j=0}^{2m}(-1)^{m+j}{2m\choose
j}\,\frac{(n-j)!(n-2m+j)!}{(n-m)!(n-2m)!}
f^{(j)}(x)f^{(2m-j)}(x)\geq 0,\qquad x\in \mathbb{R}\,.
$$
It was shown in \cite{NU(2004)} that these inequalities provide a
refinement of the Jensen inequalities for functions from the
Laguerre-P\'{o}lya class, specialized to the subclass of real-root
polynomials. In \cite{DimNik(2010)}, \eqref{e14} was deduced from
the above inequalities in the special case $m=2$, and then applied
for the derivation of certain bounds for the zeros of classical
orthogonal polynomials.

Let us substitute in \eqref{e14} $f=P_n^{(\la)}$ and
$x_0=x_{nn}(\la)$. We make use of $f(x_0)=0$ and the second order
differential equations for $f$ and $f^{\prime}$,
\begin{eqnarray*}
&& (1-x^2)\,f^{\prime\prime}-(2\la+1)x\,f^{\prime}(x)+n(n+2\la)\,f=0\,,\\
&& (1-x^2)f^{\prime\prime\prime}-(2\la+3)x\,f^{\prime\prime}(x)
+(n-1)(n+2\la+1)\,f^{\prime}=0\,,
\end{eqnarray*}
to express $f^{\prime}(x_0)$ and $f^{\prime\prime\prime}(x_0)$ in
terms of $f^{\prime\prime}(x_0)$ as follows:
\begin{eqnarray*}
&&f^{\prime}(x_0)=\frac{1-x_0^2}{(2\la+1)x_0}\,f^{\prime\prime}(x_0)\,,\\
&&f^{\prime\prime\prime}(x_0)=\Big[\frac{2\la+3)x_0}{1-x_0^2}
-\frac{(n-1)(n+2\la+1)}{(2\la+1)x_0}\Big]\,f^{\prime\prime}(x_0)\,.
\end{eqnarray*}
Putting these expressions in \eqref{e14}, canceling out the positive
factor $\big[f^{\prime\prime}(x_0)\big]^2$ and solving the resulting
inequality with respect to $x_0^2$, we arrive at the condition
$$
x_0^2\leq \frac{(n-1)(n+2\la+1)}{(n+\la)^2+3\la+\frac{5}{4}
+3\frac{(\la+1/2)^2}{n-1}}\,.
$$
Hence,
$$
x_0^2<\frac{(n-1)(n+2\la+1)}{(n+\la)^2+3\la+\frac{5}{4}} =1-
\frac{(2\la+1)(2\la+9)}{4n(n+2\la)+(2\la+1)(2\la+5)}\,.
$$
This accomplishes the proof of Theorem~C.

\section{Remarks} \label{s4}
\textbf{1.} As was mentioned in the introduction, at least for large
$n$, the bounds given in Theorem~\ref{t1}, Theorem~\ref{t3} and
Corollary~\ref{c1} are sharper than those in Theorem~A, Corollary~A
and Theorem~B, respectively. For instance, for fixed $\al,\,\be>-1$
the upper bounds for $\,1-x_{nn}(\al,\be)\,$ in Theorem~A and
Theorem~\ref{t1} are respectively
$$
\frac{(\al+1)(\al+3)}{n^2}+o(n^{-2}),\qquad
\frac{4(\al+1)(\al+2)(\al+4)}{(5\al+11)n^2}+o(n^{-2})\,,\qquad
n\to\infty\,,
$$
and
$$
(\al+1)(\al+3)-\frac{4(\al+1)(\al+2)(\al+4)}{5\al+11}=
\frac{(\al+1)^3}{5\al+11}>0\,,\qquad \al>-1\,.
$$
The same conclusion is drawn for the other two pairs of bounds when
$\la>-1/2$ is fixed and $n$ is large (it follows from the above
consideration with $\la=\al-1/2$).
\medskip

\noindent\textbf{2.} Theorem~\ref{t1} is deduced from the second
inequality in \eqref{e8} with $k=3$. Note that \eqref{e8} with $k=2$
together with \eqref{e9} and \eqref{e10} implies the estimate
$$
1-x_{nn}(\al,\be)<\frac{2(\al+1)(\al+3)}
{2n(n+\al+\be+1)+(\al+1)(\be+1)}\,,
$$
which however is less precise than the estimate in Theorem~A, and
also than the estimate of Driver and Jordaan from
\cite{DriJor(2012)},
$$
1-x_{nn}(\al,\be)<\frac{2(\al+1)(\al+3)}
{2n(n+\al+\be+1)+(\al+1)(\al+\be+2)}\,.
$$

Of course, having found the power sums $p_i(P)$, $1\leq i\leq 4$,
one could apply Proposition~\ref{p1} for derivation of lower bounds
for $1-x_{n,n}(\al,\be)$ as well. For instance, the first inequality
in \eqref{e8} with $k=4$ yields
$$
1-x_{nn}(\al,\be)>\frac{2}{\big[p_4(P)\big]^{1/4}}
$$
with $p_4(P)$ given by \eqref{e11} and $a=\al+1$, $b=\be+1$,
$t=n(n+\al+\be+1)$. However, the expression on the right-hand side
looks rather complicated to be of any use.\medskip

\noindent \textbf{3.} In \cite{GupMul(2007)} Gupta and Muldoon
proved the following upper bound for the smallest zero $x_{1n}(\al)$
of the $n$-th Laguerre polynomial $L_n^{(\al)}$:
\begin{equation}\label{e15}
x_{1n}(\al)<\frac{(\al+1)(\al+2)(\al+4)(2n+\al+1)}
{(5\al+11)n(n+\al+1)+(\al+1)^2(\al+2)}\,.
\end{equation}
Let us demonstrate how this result can be deduced from the proof of
Theorem~\ref{t1} and the well-known limit relation
$$
x_{1n}(\al)=\lim_{\be\to\infty}\frac{\be}{2}\,
\big(1-x_{nn}(\al,\be)\big)\,.
$$
Since
$$
\frac{1}{2}\big(1-x_{nn}(\al,\be)\big)\leq \frac{p_3(P)}{p_4(P)}
=\frac{a(a+1)(a+3)q_2(t)}{q_3(t)}
$$
with $a=\al+1$, $b=\be+1$, $t=n(n+\al+\be+1)$, and $q_2(t)$,
$q_3(t)$ given in \eqref{e10} - \eqref{e11}, we have
\begin{equation}\label{e16}
x_{1n}(\al)=\lim_{\be\to\infty}\frac{\be}{2}\,
\big(1-x_{nn}(\al,\be)\big)\leq a(a+1)(a+3)
\lim_{b\to\infty}\,\frac{b\,q_2(t)}{q_3(t)}\,.
\end{equation}
Using
$$
\lim_{b\to\infty}\,\frac{t}{b}=n
$$
and the explicit form of $q_2(t)$ and $q_3(t)$, we find
\[
\begin{split}
\lim_{b\to\infty}\,\frac{b\,q_2(t)}{q_3(t)}&=
\lim_{b\to\infty}\,\frac{q_2(t)/b^2} {q_3(t)/b^3}\\
&=\frac{2n^2+3a\,n+a^2}{(5a+6)n^3+2a(5a+6)n^2+a^2(6a+7)n+a^3(a+1)}\\
&=\frac{2n+a+1}{(5a+6)n(n+a)+a^2(a+1)}\,.
\end{split}
\]
By substituting the latter expression in \eqref{e16} and setting
$a=\al+1$, we arrive at \eqref{e15}\,.
\medskip

\noindent \textbf{4.} We already mentioned in the introduction that,
for every fixed $\la>-1/2$, the ratio $\,r(\la,n)\,$ of the upper
and the lower bound for $\,1-x_{nn}^2(\la)$, given by
Theorem~\ref{t1} and Theorem~C, respectively, tends to a limit which
does not exceed $1.6$. More precisely,
$$
r(\la,n)=\varrho(\la)\psi(\la,n)\,,
$$
where
$$
\varrho(\la)=\frac{8(2\la+3)(2\la+7)}{(2\la+9)(10\la+17)}\,,\qquad
\varphi(\la,n)=\frac{n(n+2\la)+(2\la+1)(2\la+5)/4}{n(n+2\la)+(2\la+1)^2/8}\,.
$$

\begin{figure}[h]
\centering
\resizebox*{6cm}{!}{\includegraphics{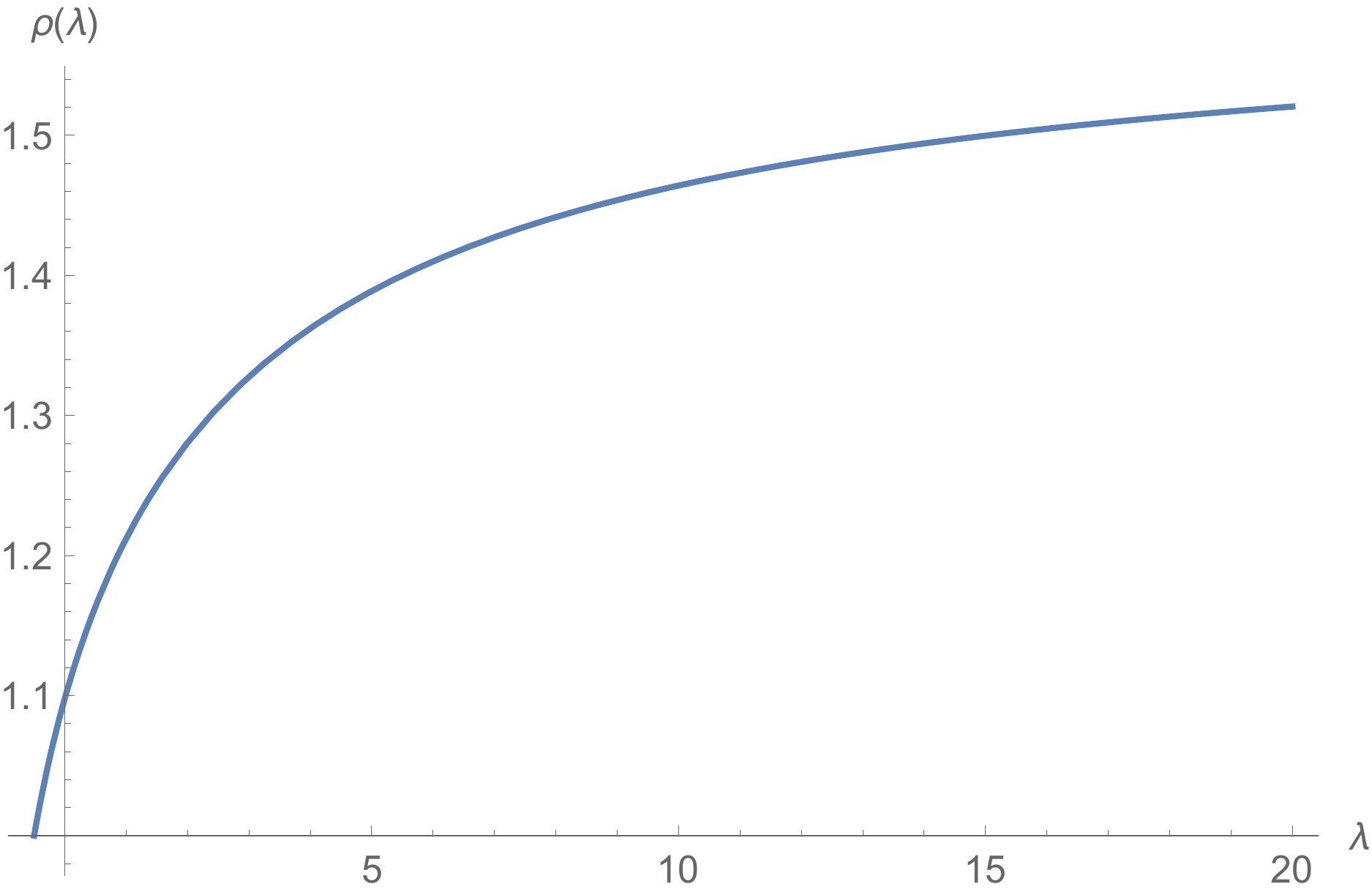}}\hspace{5pt}
\caption{The graph of $\varrho(\la)$.} \label{sample-figure}
\end{figure}

The function $\varrho(\la)$ is monotonically increasing in the
interval $(-1/2,\infty)$ assuming values between $1$ and $1.6$ (see
Fig. 1) while, for a fixed $\la>-1/2$,
$\lim_{n\to\infty}\varphi(\la,n)=1$.

\medskip

\noindent \textbf{5.} The Euler--Rayleigh approach assisted with
symbolic algebra has been applied in \cite{NU(2019)} to the
derivation of bounds for the extreme zeros of the Laguerre
polynomials, and in \cite{AN2018, NS(2017), NU(2017)} to the
estimation of the extreme zeros of some non-classical orthogonal
polynomials, which are related to the sharp constants in some
Markov-type inequalities in weighted $L_2$ norms.

\section*{Acknowledgements} The author would like to thank the Isaac
Newton Institute for Mathematical Sciences, Cambridge, for support
and hospitality during the programme \emph{Approximation, sampling
and compression in data science} where work on this paper was
undertaken. This work was partially supported by EPSRC grant no
EP/K032208/1, by a grant from the Simons Foundation, and by Sofia
University Research Fund under Contract 80-10-17/2019.

\end{document}